\documentclass{birkmult}
\usepackage[cp1251]{inputenc}
\usepackage{amsmath,amsthm,amssymb,graphicx,color}
\usepackage[colorlinks,pdfauthor={Yu. Brezhnev},
            pdfwindowui=false,
            bookmarksopen=false,bookmarks=false]{hyperref}

\def~{\unskip\nobreak\hspace{0.27778em}\ignorespaces}
\def\,{\ifmmode\mskip+1.5mu\else\kern+.08333em\fi\relax}
\def\!{\ifmmode\mskip-1.5mu\else\kern-.08333em\fi\relax}
\newdimen{\FontSize} \delimiterfactor=990
\def\cdot{\,{\mathchar"2201}\,}

\makeatletter
\def\dot#1{\mathbf{\mathaccent"705F\mathnormal{#1}}}
\def\ddot#1{\mathbf{\mathaccent"707F\mathnormal{#1}}}
\def\dddot#1{{\mathop{{}#1}\limits^{\vbox to-1.5\ex@{\kern-\tw@\ex@
\hbox{\larger[2]\rm.\kern-0.12em.\kern-0.12em.}\vss}}}} \makeatother

\def\ds{\displaystyle}

\def\s{\scriptscriptstyle}
\def\ie{i.\,\,e.}

\def\eg{\text{e.\,\,g.}}
\def\re{{\,\mathrm{e}}}
\def\ri{{\mathrm{i}}}

\def\wpp {\wp\Smaller{'}}

\def\II{{\textbf{\textsc{I\kern -0.15em I}}}}
\def\III{{\textbf{\textsc{I\kern -0.17em I\kern -0.17em I}}}}
\def\IIIab{\III\!\raise0.2ex\hbox{$\begin{smallmatrix}
\alpha\\[-0.1ex]\varrho\end{smallmatrix}$}\kern-0.1em}

\def\DEF{\mathrel{\vcenter{\hbox{$:$}}{=}}}

\def\=={\mathrel{\phantom{=}}}

\newcommand{\Dtheta   }[1][1]{{\ds\theta{\vbox to 1.7ex{}}^\prime
\nobreak{}_{\mskip-7.5mu#1}}\mskip-1.5mu{\vbox to 1.7ex{}}\relax}

\def\?{\textrm{\protect\footnotesize$\red\mathchar"446$}}

\def\vphi{\smash[b]{\raise-0.22ex\hbox{\Smaller[2]{\s\boldsymbol-}}
  \mkern-8.00mu\raise0.435ex\hbox{\scalebox{1}[0.9]{$\varphi$}}}\relax}
\def\vpsi{\smash[b]{\raise-0.22ex\hbox{\Smaller[2]{\s\boldsymbol-}}
  \mkern-8.95mu\raise0.408ex\hbox{\scalebox{1}[0.9]{$\psi   $}}}\relax}

\def\END{\end{document}}

\def\Mfrac#1#2{\def\FRAC{\hbox{\smaller[1]$\ds\frac{#1}{#2}$}}
 \ifdim\FontSize=10pt\raise0.057ex\FRAC\else
{\ifdim\FontSize=11pt\raise0.050ex\FRAC\else\raise0.051ex\FRAC\fi}\fi}

\def\mfrac#1#2{\def\FRAC{\hbox{\smaller[2]$\ds\frac{#1}{#2}$}}
 \ifdim\FontSize=10pt\raise0.1166ex\FRAC\else
{\ifdim\FontSize=11pt\raise0.104ex\FRAC\else\raise0.0968ex\FRAC\fi}\fi}

\newcommand{\AB}[2]{\nobreak\mskip-1mu
\raise0.12ex\hbox{\smaller[2]\renewcommand{\arraystretch}{0.42}%
\setlength{\arraycolsep}{0pt}\raise0.16ex\hbox{$\big[$}%
$\begin{array}{c}\scriptstyle#1\\ \scriptstyle#2\end{array}$%
\raise0.16ex\hbox{$\big]$}}\mskip-1mu\relax}

\newcommand{\mbig}[2][4]{\vcenter{\hbox{%
\ifcase#1$#2$\or\small$\big#2$\or$\big#2$\or\large$\big#2$%
\or\footnotesize$\Big#2$\or\small$\Big#2$\or$\Big#2$\or\large$\Big#2$%
\or$\bigg#2$\or\large$\bigg#2$\or$\Bigg#2$\or\large$\Bigg#2$%
\or\Large$\Bigg#2$\else\huge$\bigg#2$\fi}}\relax}

\newdimen{\upp} \newdimen{\Dim} \newdimen{\DDim}
\newcommand{\pow}[4][\displaystyle]{
\settowidth{\upp}{\hbox{$#1 {}^{#4}$}}%
\settowidth{\Dim}{\hbox{$#1{#2}$}}%
\settowidth{\DDim}{\hbox{$#1{#2}_{#3}$}}%
\addtolength{\DDim}{-\Dim} {#1{#2}_{#3}\kern-\DDim {}^{#4}}
\addtolength{\DDim}{-\upp}\ifdim\DDim>0pt \kern+\DDim \fi\relax}

\makeatletter 
\providecommand{\larger}[1][1]{
\ifmmode\relax\else\@tempcnta\ifx\@currsize\normalsize 4\else
 \ifx\@currsize\small 3\else\ifx\@currsize\footnotesize 2\else
 \ifx\@currsize\large 5\else\ifx\@currsize\Large 6\else
 \ifx\@currsize\LARGE 7\else\ifx\@currsize\scriptsize 1\else
 \ifx\@currsize\tiny  0\else\ifx\@currsize\huge 8\else
 \ifx\@currsize\Huge  9\else 4 \fi\fi\fi\fi\fi\fi\fi\fi\fi\fi
\advance\@tempcnta#1\relax\ifnum\@tempcnta<\z@ \@tempcnta\z@ \fi
\ifcase\@tempcnta\tiny\or\scriptsize\or\footnotesize\or\small
 \or\normalsize\or\large\or\Large\or\LARGE\or\huge\else\Huge\fi\fi}
\makeatother

\providecommand{\smaller}[1][1]{\larger[-#1]}

\newcommand{\Smaller}[2][1]{\hbox{\smaller[#1]$#2$}}

\makeatletter \def\@dotsep{.8} \makeatother 

\def\red{\color{red}}
\def\blue{\color{blue}}

\newcommand{\Width}[2][\displaystyle]
{\settowidth\Dim{\hbox{$#1#2$}}\texttt{\tiny\blue$($\the\Dim$)$}}
\newcommand{\Height}[2][\displaystyle]
{\settoheight\Dim{\hbox{$#1#2$}}\texttt{\tiny\blue$($\the\Dim$)$}}
\newcommand{\Depth}[2][\displaystyle]
{\settodepth\Dim{\hbox{$#1#2$}}\texttt{\tiny\blue$($\the\Dim$)$}}
\def\Psix{\mbox{$\mathcal{P}_6$}}
\def\u{{\mathfrak{u}}}
\def\Q{{\mathcal{Q}}}
\def\R{\mathcal{R}}

\author[Yu.~Brezhnev]{Yurii V.~Brezhnev}

\title[The Painlev\'e transcendents and uniformization]
{The sixth Painlev\'e transcendent\\
and uniformizable orbifolds}

\thanks{Research supported by the Federal Targeted Program under
contract 02.740.11.0238.}

\begin{document}
\hfill {\smaller[2]{\bf Painlev\'e equations and Related Topics,}
193--198. {\smaller \copyright} De~Gruyter 2012}\bigskip

\maketitle

\vspace{-1.7em}
\section{Algebraic solutions of \Psix\ and uniformization theory}
The sixth Painlev\'e transcendent
\begin{equation*}
\begin{aligned}
\Psix\!:\;\; y_{\mathit{xx}}^{}={}&\frac12\!
\left(\frac1y+\frac{1}{y-1}+ \frac{1}{y-x}\right)\! y_x^2-
\left(\frac1x+\frac{1}{x-1}+
\frac{1}{y-x}\right)\! y_x^{}\\[0.5em]
&{+}\;\frac{y\,(y-1)(y-x)}{x^2(x-1)^2}\!
\left\{\alpha-\beta\,\frac{x}{y^2}+\gamma\,\frac{x-1}{(y-1)^2}-
\Big(\delta-\mfrac12\Big)\,\frac{x\,(x-1)}{(y-x)^2} \right\}
\end{aligned}
\end{equation*}
is known to be a rich source of nontrivial algebraic solutions $y=f(x)$
and genera of these solutions, as genera of corresponding algebraic
curves $F(x,y)=0$, may be made as great as is wished. The relation of
such solutions to the uniformization theory is based on the
$\wp$-representation of the \Psix:
\begin{equation}\label{P6wp}
-\frac{\pi^2}{4}\,\frac{d^2z}{d\tau^2}=
\alpha\,\wpp(z|\tau)+\beta\,\wpp(z-1|\tau)+\gamma\,\wpp(z-\tau|\tau)+
\delta\,\wpp(z-1-\tau|\tau)
\end{equation}
obtainable via the transcendental change $(x,y)\mapsto(z,\tau)$
(Painlev\'e (1906), Manin--Babich--Bordag (1996)):
\begin{equation}\label{subswp}
x=\frac{\vartheta_4^4(\tau)}{\vartheta_3^4(\tau)}\,,\qquad
y=\frac13+\frac13\frac{\vartheta_4^4(\tau)}{\vartheta_3^4(\tau)}
-\frac{4}{\pi^2} \frac{\wp(z|\tau)}{\vartheta_3^4(\tau)}\,.
\end{equation}
Thus, knowledge of $z(\tau)$-dependence leads to a parametric
representation for solution $y=f(x)$ and, in particular, to parametric
representation of algebraic solutions. In their full generality these
dependencies are known for the Picard--Hitchin class of solutions. For
example, Picard's case $\alpha=\beta=\gamma=\delta=0$ corresponds to
$z=A\,\tau+B$. In Hitchin's case $\alpha=\beta=\gamma=\delta=\frac18$
the dependence $z(\tau)$ is more complicated (obtainable through
Okamoto's transformations) but parametric form of solution is, however,
found to be very compact
\begin{equation}\label{ph}
y_{\mbox{\tiny Pic}}^{}=-\frac{\vartheta_4^2(\tau)}
{\vartheta_3^2(\tau)}\frac{\theta_2^2}{\theta_1^2}\,,
\qquad
y_{\mbox{\tiny Hit}}^{}=
\frac{\vartheta_4^2(\tau)}{\vartheta_3^2(\tau)}
\left\{\pi\,\frac{\vartheta_2^2(\tau)\!\cdot\!
\theta_2\,\theta_3\,\theta_4}
{\Dtheta+2\,\pi A\,\theta_1}-\theta_2^2\right\}
\!\frac{1}{\theta_1^2}\,,
\end{equation}
where $\theta$'s are understood to be equal to
$\theta_k(A\,\tau+B|\tau)$ with arbitrary constants $A$, $B$ and
$\Dtheta\DEF\Dtheta(A\,\tau+B|\tau)$. Purely algebraic solutions
correspond to $A\,\tau+B=\frac{\nu}{N}\,\tau+\frac{\mu}{N}$ with
integral $\nu$, $\mu$, and $N$.

Uniformizing functions are known to be determined in terms of the
auxiliary 2nd order linear Fuchsian ODEs
$\Psi_{\mathit{yy}}=\frac12\,\Q(x,y)\,\Psi$, where $\Q$, as a rational
function of $x$ and $y$, contains all the information about
corresponding Riemann surface $\R$ (or orbifold $\mathfrak{T}$). Since
the function $x=\chi(\tau)$ in \eqref{subswp} is the very well-known
one and its Fuchsian $\Gamma(2)$-equation
$\Psi''=-\frac14\,\frac{x^2-x+1}{x^2(x-1)^2}\,\Psi$ is also known, we
obtain nontrivial (solvable) Fuchsian equations for the second
uniformizing function $y(\tau)$. Manipulations with Fuchsian equations
themselves are not convenient because we constantly handle the
multivalued functions-inversions; the ratios like
$\tau=\Psi_1(x)/\Psi_2(x)$. For this reason we invert the standard
Schwarz derivative $\{\tau,x\}$ into the `reverse' object
$[x,\tau]=-\{\tau,x\}$ and work with the autonomous  ODEs
\begin{equation}\label{Qy}
[y,\tau]=\Q(x,y)\,,\qquad\text{where\ \ }[y,\tau]\DEF
\frac{\dddot{\smash[b]{y}}}{\dot y^3}-
\frac32\,\frac{{\ddot y}^2}{{\dot y}^4}\,,
\end{equation}
defining uniformizing single-valued functions and other single-valued
objects.

\section{On the general solution to equation \eqref{P6wp}}
Complete structure of the analytic continuations (a connection problem)
of arbitrary solutions to \Psix\ is the subject matter of the series
works by D.~Guzzetti (see, \eg, \cite{guz}). Analyzing these results,
it would appear reasonable that the ramification structure of all (not
necessarily algebraic)  solutions to the \Psix-equation in the vicinity
of critical points is described by a function series of the kind
$$
y= A+R\big[(x-e)^a\ln^n \!(x-e)\big]+\cdots\,,
$$
where $e=\{0,1\}$, $a\in\mathbb{C}$, $n\in\mathbb{Z}$, and $R[\ldots]$
is a rational function of its argument. In the language of uniformizing
Painlev\'e substitution \eqref{subswp} this point is self-suggested: in
the upper $(\tau)$-half-plain $\mathbb{H}^+$ the $x$-function has an
exponential behavior in the neighborhood of the points
$x=\{0,1,\infty\}$:
$$
x\stackrel{\tau\to0}{=}0+16\,
\exp\!\left(\mfrac{\pi}{\ri\,\tau}\right)+\cdots\,,\quad
x\stackrel{\tau\to\infty}{=}
1-16\,\re^{\pi\ri\tau}+\cdots\,,\quad
x\stackrel{\tau\to1}{=}\frac{1}{16}\,
\exp\!\left(\mfrac{\pi\,\ri}{\tau-1}\right)+\cdots
$$
(the uniformizing $\tau$-parameter itself is defined up to a
fraction-linear transformation). It follows (the conjecture) that the
$y$-function has also the single-valued character about each of the
branch-point pre-images:
\begin{equation}\label{str}
y(\tau)=A+B\,(\tau-\tau_\mathrm{o})^n
\exp\!\Big(\mfrac{-a\,\pi\ri}{\tau-\tau_\mathrm{o}^{}}\Big)+\cdots\,,
\qquad
y(\tau)=A+B\,\tau^n\re^{a\pi\ri\tau}+\cdots\,.
\end{equation}
as $\tau\to \tau_\mathrm{o}\in\mathbb{R}$ or, respectively,
$\tau\to+\ri\infty$. For example, all asymptotics appearing in
\cite{guz} fit this behavior. We can therefore rewrite Eqs.~\Psix\ and
\eqref{P6wp} in form of modification of purely `algebraic' uniformizing
Schwarz--Fuchs 3rd order ODE \eqref{Qy}:
\begin{equation}\label{conj}
[y,\tau]=A\,y_x^{\s4}+B\,y_x^{\s3}+C\,y_x^{\s2}+
D\,y_x^{\s1}+E\,,
\end{equation}
where $(A,B,C,D,E)$ are certain rational functions of $x,y$ and
quadratic polynomials in parameters $(\alpha,\beta,\gamma,\delta)$
(explicit expressions are too cumbersome to display here). Because of
outstanding character of \Psix, this equation may be treated as a
generator of `infinite genus curves'. In the case of algebraic
solutions the right hand side of Eq.\;\eqref{conj} becomes a rational
function $\Q(x,y)$, that is \eqref{Qy}. We conjecture that all the
Painlev\'e solutions to Eq.\;\eqref{conj} are the globally
single-valued analytic functions with the structure \eqref{str} and the
domain of their existence is a half-plain (under suitable normalization
of $\tau$). It is known that solutions to the lower Painlev\'e
equations $\mathcal{P}_{1\ldots5}$ (under an appropriate modification
\cite{gromak}) are the single-valued functions on
$\overline{\mathbb{C}}$. In this respect, the pass from \Psix-equation
over $\overline{\mathbb{C}}\backslash\{0,1,\infty\}$ to the
$\mathbb{H}^+$ and uniformization theory related to the coverings of a
three punctured $\Gamma(2)$-orbifold becomes very natural.

\vspace{-1em}
\section{Calculus: Abelian integrals and affine (analytic)
connections}

Insomuch as we have not only $\tau$-representations for the scalar
(\ie\ automorphic) functions on  $\R$'s but rules for differential
computations with theta-functions of arbitrary arguments \cite{br}  we
can close the differential apparatus on orbifolds $\mathfrak{T}$ whose
compactifications are corresponding Painlev\'e $\R$'s. This includes
the additively automorphic functions (Abelian integrals),
differentials, and  covariant differentiation, say, of 1-differentials
$\nabla=\partial_\tau-\Gamma(\tau)$. The latter leads to necessity to
introduce the geometric connection object $\Gamma(\tau)$, which
transforms according to the standard rule
$\widetilde{\Gamma}(\tilde\tau)\,d\tilde\tau=
\Gamma(\tau)\,d\tau-d\ln\!\frac{d\tilde\tau}{d\tau}$ under
$\mathrm{SL}_2(\mathbb{R})$-transformations and respects the factor
topology of $\mathbb{H}^+\!\!\big/\pi_1(\mathfrak{T})$. The
characteristic feature of the (complex) 1-dimensional case (orbifolds
and Riemann surfaces) is that it is completely described by the
invariant 3rd order ODE \eqref{Qy}. Therefore closed collection of data
for the theory is given by the set $\big\{y(\tau),\,\dot
y(\tau),\,\ddot y(\tau)\big\}$ if, however, the automorphism group of
the generator $y(\tau)$ coincides with $\pi_1(\mathfrak{T})$. In
general, automorphisms of the field generators are not bound to
coincide with $\pi_1(\mathfrak{T})$ since the choice of the pair
$(x,\,y)$ is not unique. It is found however that the set of Painlev\'e
orbifolds coming from Picard--Hitchin's curves \eqref{ph} is not the
case: $\mathrm{Aut}\,y(\tau) \cong \pi_1(\mathfrak{T})$. In this regard
the many Painlev\'e curves (we suggest that all) stand out majority of
classical modular equations originated from purely group-algebraic
considerations related to the group $\mathrm{PSL}_2(\mathbb{Z})$ or
some its subgroups. By this means the expression
$$
\Gamma(\tau)=\frac{d}{d\tau}\!\ln\dot y(\tau)+
\text{arbitrary (Abelian) 1-differential}
$$
provides a general form of the sought-for connection on Painlev\'e
$\mathfrak{T}$. We can normalize this $\Gamma(\tau)$ to have only first
order poles (residues) and, integrating the transformation law above,
one can see that the sum of such residues is invariant
$$
\int_{\partial\R}\!\!\!\!\widetilde{\Gamma}
(\tilde\tau)\,d\tilde\tau
=\int_{\partial\R}\!\!\!\!\Gamma(\tau)\,d\tau=
(2\,g-2)\cdot2\,\pi\,\ri\,;
$$
it depends only on genus and, in effect, is equal to the number of
zeroes of a holomorphic differential $\dot\u(\tau)$. Varying the
holomorphic differentials $\dot{\u}_k(\tau)$ we can impart the simpler
from to the connection
$$
\Gamma(\tau)=\frac{d}{d\tau}\!\ln\dot{\u}(\tau)+
\sum\limits_{k=1}^{g}\dot{\u}_k(\tau)
$$
and build the elementary $\Gamma$ with a single pole (if genus $g>1$
then the analytic connection does always have a singularity). So we
have  the set of invariant objects $\big\{y(\tau),\,\dot
y(\tau),\,\Gamma(\tau)\big\}$ since functions $x(\tau)$, $y(\tau)$ are
completely at hand. The remarkable fact is that \textit{affine
$($analyitc\/$)$ connection on an arbitrary $\mathfrak{T}$ satisfies an
autonomous  ODE\/ $\Xi(\dddot{\smash[b]{\,\Gamma}},\ddot{\Gamma},
\dot{\Gamma},\Gamma)=0$} and there is an algorithm how to derive it.

For completeness we should involve into analysis the integrals of
closed \mbox{1-forms} on our $\mathfrak{T}$'s and $\R$'s, if only
because there are exact 1-forms whose integrals lead to the scalar
objects. On the other hand, uniformization of any higher genera curves
is reduced to uniformization of \textit{zero} genus orbifolds and the
latter form towers and hierarchies. In the Painlev\'e uniformizing
theory, in one way or another, many classical and nonclassical zero
genus known orbifolds appear \cite{br}. In turn they are related to
nonzero genus curves which may cover elliptic ones, \ie\ tori. We thus
obtain a possibility to construct explicitly Abelian integrals if they
come from an elliptic cover. Here is a good example along these lines.

The Chudnovsky orbifold defined by the Fuchsian equation
$(z^3-z)\,\Psi''+(3\,z^2-1)\,\Psi'+z\,\Psi=0$ is related, through the
Halphen transformation (zero genus elliptic cover) $z=\wp(\u)$, to the
Fuchsian equation on the lemniscatic torus $\wpp{}^2=4\,\wp^3-4\,\wp$.
Correlating these facts we derive the nice $\tau$-representation for
the everywhere finite object $\u$ and analog of \eqref{Qy}---the
uniformizing Schwarz equation:
$$
[\u,\tau]=-2\,\wp(2\,\u)\,,\qquad \u(\tau)=
\frac12\frac{\vartheta_3(\tau)}{\vartheta_2(\tau)}\!\cdot\!
{}_2F_1\!\!\!\left(\Mfrac12,\Mfrac14;\Mfrac54\Big|
\Mfrac{\vartheta_3^4(\tau)}{\vartheta_2^4(\tau)}\right)
$$
(the check is a good exercise). This is a first \textit{explicit and
analytic} $\tau$-representation for an additively automorphic function
(Abelian integral $\u=\wp^{\s1}\!(z)$) on an orbifold (Riemann surface)
of a \textit{negative} curvature $-1$. Under suitable cover this
$u(\tau)$ may produce the $\tau$-representation for $\u$-integrals on
higher genus curves; examples of the analogous ODEs and their solutions
can also be obtained. All of them can be related to the Painlev\'e
curves.

\end{document}